\documentclass[12pt]{amsart}

\usepackage{amssymb,amsmath}
\usepackage{latexsym}
\usepackage{pstricks}
\usepackage{rotate}
\newtheorem{example}{Example}[section]
\newtheorem{note}[example]{Note}

\newtheorem{theorem}[example]{Theorem}

\newtheorem{definition}[example]{Definition}
\newtheorem{proposition}[example]{Proposition}

\newtheorem{lemma}[example]{Lemma}
 
\def\Proof{\noindent \it Proof -- \rm}                                                 
\def\qed{\hspace{3.5mm} \hfill \vbox{\hrule height 3pt depth 2 pt width 2mm}
\bigskip}
\def\PW{{\rm PW}}
\def\std{{\rm std}}
\def\q{{\bf q}}
\def\t{{\bf t}}
\def\Sym{{\bf Sym}}
\def\QSym{{\it QSym}}
\def\NCSF{{\bf Sym}}
\def\FSym{{\bf FSym}}
\def\maj{{\rm maj}}
\def\bo{{\small \,\Box}}

\def\<{\langle}
\def\>{\rangle}

\def\SG{{\mathfrak S}}

\def\t{{\bf t}}

\def\tHH{{\rm \tilde H}}

\def\tGG{{\rm \tilde G}}

\def\Des{{\rm Des\,}}
\def\bkt{{\bf \tilde k}}

\def\KK{{\mathcal K}}

\def\tHH{{\rm \tilde H}}
\def\btg{{\bf \tilde g}}
\def\y{{\bf y}}

\def\ev{{\rm ev}}

\newdimen\squaresize
\newdimen\thickness        
                                                    
\def\square#1{\hbox{\vrule width \thickness
   \vbox to \squaresize{\hrule height \thickness\vss                            
      \hbox to \squaresize{\hss#1\hss}
   \vss\hrule height\thickness} 
\unskip\vrule width \thickness} 
\kern-\thickness}                                                            
                               
\def\vsquare#1{\vbox{\square{$\casestyle#1$}}\kern-\thickness}
\def\blank{\omit\hskip\squaresize}

\def\young#1{\vcenter{%
    \squaresize=18pt\thickness=0.5pt\let\casestyle=\relax
    \vbox{\smallskip\offinterlineskip
      \halign{&\vsquare{##}\cr #1}}}}

\def\moyyoung#1{\vcenter{%
    \squaresize=28pt\thickness=0.6pt\let\casestyle=\relax
    \vbox{\smallskip\offinterlineskip
      \halign{&\vsquare{##}\cr #1}}}}

\def\bigyoung#1{\vcenter{%
    \squaresize=38pt\thickness=0.8pt\let\casestyle=\displaystyle
    \vbox{\smallskip\offinterlineskip
      \halign{&\vsquare{##}\cr #1}}}}

\def\boxit#1#2{\setbox1=\hbox{\kern#1{#2}\kern#1}%
\dimen1=\ht1 \advance\dimen1 by #1 \dimen2=\dp1 \advance\dimen2 by #1
\setbox1=\hbox{\vrule height\dimen1 depth\dimen2\box1\vrule}%
\setbox1=\vbox{\hrule\box1\hrule}%
\advance\dimen1 by .4pt \ht1=\dimen1
\advance\dimen2 by .4pt \dp1=\dimen2 \box1\relax}
\def\bo#1{\boxit{1pt}{$#1$}}

\def\rdet{{\rm rdet\,}}

\title{Noncommutative symmetric functions with matrix parameters}

\author[A. Lascoux, J.-C. Novelli, J.-Y. Thibon]%
{Alain Lascoux, Jean-Christophe Novelli, and Jean-Yves Thibon}

\date{\today}

\begin{document}

\begin{abstract}
We define new families of noncommutative symmetric functions and
quasi-symmetric functions depending on two matrices of parameters, and more
generally on parameters associated with paths in a binary tree.
Appropriate specializations of both matrices then give back the two-vector
families of Hivert, Lascoux, and Thibon %[arXiv: math.CO/0106191]
and the noncommutative Macdonald functions of Bergeron and Zabrocki.
%[Discr. Math. {\bf 298} (2005), 79--103].
\end{abstract}

\maketitle

%%%%%%%%%%%%%%%%%%%%%%%%%%%%%%%%%%%%%%%%%%%%%%%%%%%%%%%%%%%%%%%%%%%%%%%%%%%%%%
%%%%%%%%%%%%%%%%%%%%%%%%%%%%%%%%%%%%%%%%%%%%%%%%%%%%%%%%%%%%%%%%%%%%%%%%%%%%%%
%%%%%%%%%%%%%%%%%%%%%%%%%%%%%%%%%%%%%%%%%%%%%%%%%%%%%%%%%%%%%%%%%%%%%%%%%%%%%%
\section{Introduction}

The theory of Hall-Littlewood, Jack, and Macdonald polynomials is one of the
most interesting subjects in the modern theory of symmetric functions. 
It is well-known that combinatorial properties of symmetric functions can be
explained by lifting them to larger algebras (the so-called combinatorial Hopf
algebras), the simplest examples being $\NCSF$ (Noncommutative symmetric
functions \cite{NCSF1}) and its dual $QSym$ (Quasi-symmetric functions
\cite{Ges}).

There have been several attempts to lift Hall-Litttlewood and Macdonald
polynomials to $\NCSF$ and $QSym$ \cite{BZ,Hiv,HLT,NTW,Tev}. The analogues
defined in \cite{BZ} were similar to,  though different from, those of
\cite{HLT}.  These last ones admitted multiple parameters $q_i$ and $t_i$,
which however could not be specialized to recover the version of \cite{BZ}.

The aim of this article is to show that many more parameters can be introduced
in the definition of such bases. Actually, one can have a pair of $n\times n$
matrices $(Q_n,T_n)$ for each degree $n$. The main properties established in
\cite{BZ} and \cite{HLT} remain true in this general context, and one recovers
the BZ and HLT polynomials for appropriate specializations of the matrices.

In the last section, another possibility involving quasideterminants is
explored.
One can then define bases involving two almost-triangular matrices of
parameters in each degree. We shall see on some examples that if these
matrices are chosen such as to give a special basis for the row and column
compositions, special properties arise for hook compositions $(n-k,1^k)$. For
example, on can obtain a basis whose commutative image reduces the the
Macdonald $P$-polynomials for hook compositions.

One should not expect that constructions at the level of $\Sym$ and $QSym$
could lead to general results on ordinary Macdonald polynomials.  Even for
Schur functions, one has to work in the algebra of standard tableaux, $\FSym$,
to understand the Littlewood-Richardson rule. However, the \emph{analogues} of
Macdonald polynomials which can be defined in $\Sym$ and $QSym$ have
sufficiently much in common with the ordinary ones so as to suggest
interesting ideas.  The most startling one is that the usual Macdonald
polynomials could be specializations of a family of symmetric functions with
many more parameters%
\footnote{This idea has been explored in an unpublished work~\cite{HLTmulti},
in which rather convincing $({\bf q}, {\bf t})$-Kostka matrices have been
constructed up to $n=5$.}.

%%%%%%%%%%%%%%%%%%%%%%%%%%%%%%%%%%%%%%%%%%%%%%%%%%%%%%%%%%%%%%%%%%%%%%%%%%%%%%
%%%%%%%%%%%%%%%%%%%%%%%%%%%%%%%%%%%%%%%%%%%%%%%%%%%%%%%%%%%%%%%%%%%%%%%%%%%%%%
%%%%%%%%%%%%%%%%%%%%%%%%%%%%%%%%%%%%%%%%%%%%%%%%%%%%%%%%%%%%%%%%%%%%%%%%%%%%%%
\section{Notations}

Our notations for noncommutative symmetric functions will be as
in~\cite{NCSF1,NCSF2}. 
Here is a brief reminder.

The Hopf algebra of noncommutative symmetric
functions is denoted by $\Sym$, or by $\Sym(A)$ if we consider the realization
in terms of an auxiliary alphabet. Bases of $\Sym_n$ are labelled by
compositions $I$ of $n$. The noncommutative complete and elementary functions
are denoted by $S_n$ and $\Lambda_n$, and the notation $S^I$ means
$S_{i_1}\dots S_{i_r}$. The ribbon basis is denoted by $R_I$.
The notation $I\vDash n$ means that $I$ is a composition of $n$.
The conjugate composition is denoted by $I^\sim$.

The graded dual of $\Sym$ is $QSym$ (quasi-symmetric functions).
The dual basis of $(S^I)$ is $(M_I)$ (monomial), and that of $(R_I)$
is $(F_I)$. The \emph{descent set} of $I=(i_1,\dots,i_r)$ is
$\Des(I) = \{ i_1,\ i_1+i_2, \ldots , i_1+\dots+i_{r-1}\}$.

%%%%%%%%%%%%%%%%%%%%%%%%%%%%%%%%%%%%%%%%%%%%%%%%%%%%%%%%%%%%%%%%%%%%%%%%%%%%%%
%%%%%%%%%%%%%%%%%%%%%%%%%%%%%%%%%%%%%%%%%%%%%%%%%%%%%%%%%%%%%%%%%%%%%%%%%%%%%%
%%%%%%%%%%%%%%%%%%%%%%%%%%%%%%%%%%%%%%%%%%%%%%%%%%%%%%%%%%%%%%%%%%%%%%%%%%%%%%
\section{$\Sym_n$ as a Grassmann algebra}

Since for $n>0$, $\Sym_n$ has dimension $2^{n-1}$, it can be identified
(as a vector space) with a Grassmann algebra on $n-1$ generators
$\eta_1,\dots,\eta_{n-1}$
(that is, $\eta_i\eta_j=-\eta_j\eta_i$). 
This identification is meaningful, for example, in the context of the
representation theory of the $0$-Hecke algebras $H_n(0)$.
Indeed (see~\cite{NCSF6}), the quiver of $H_n(0)$ admits a simple description
in terms of this identification. 

If $I$ is a composition of $n$ with descent set $D=\{d_1,\dots,d_k\}$,
we make the identification
\begin{equation}
R_I \longleftrightarrow\ \eta_D:=\eta_{d_1}\eta_{d_2}\dots \eta_{d_k}\,.
\end{equation}
For example, $R_{213}\leftrightarrow \eta_2\eta_3$.
We then have
\begin{equation}
S^I \longleftrightarrow\ (1+\eta_{d_1})(1+\eta_{d_2})\dots(1+\eta_{d_k})
\end{equation}
and
\begin{equation}
\Lambda^I \longleftrightarrow\ \prod_{i=1}^{n-1} \theta_i,
\end{equation}
where $\theta_i=\eta_i$ if $i\not\in D$ and $\theta_i=1+\eta_i$ otherwise.
Other bases have simple expression under this identification, \emph{e.g.},
\begin{equation}
\Psi_n \longleftrightarrow
  1-\eta_1+\eta_1\eta_2 - \dots + (-1)^{n-1} \eta_1\dots\eta_{n-1},
\end{equation}
and
\begin{equation}
\Phi_n \longleftrightarrow \sum_{k=0}^{n-1} \frac{(-1)^k}{\binom{n-1}{k}} E_k,
\end{equation}
where $E_k = \sum_{j_1<\dots<j_k} \eta_{j_1} \dots \eta_{j_k}$.
The $q$-Klyachko element \cite{NCSF1}
\begin{equation}
K_n(q) = \sum_{I\vdash n} q^{\maj(I)} R_I,
\end{equation}
is
\begin{equation}
(1+q\eta_1)(1+q^2\eta_2) \dots (1+q^{n-1}\eta_{n-1}),
\end{equation}
and Hivert's Hall-Littlewood basis \cite{Hiv} is
\begin{equation}
H_I(q) := 
(\eta_{d_1}+q)(\eta_{d_2}+q^2) \dots (\eta_{q{_n-1}}+q^{n-1}).
\end{equation}

%%%%%%%%%%%%%%%%%%%%%%%%%%%%%%%%%%%%%%%%%%%%%%%%%%%%%%%%%%%%%%%%%%%%%%%%%%%%%%
\subsection{Structure on the Grassmann algebra}

Let $*$ be the anti-involution given by $\eta_i^*=(-1)^i\eta_i$.
Recall that the Grassmann integral is defined by 
\begin{equation}
\int d\eta\, f :=
   f_{12\dots n-1},\qquad\text{where}\quad
 f=\sum_{k}\sum_{i_1<\dots<i_k}f^{i_1\dots i_k}\eta_{i_1}\dots\eta_{i_k}.
\end{equation}
We define a bilinear form on $\Sym_n$ by
\begin{equation}
(f,g)=\int d\eta\, f^* g
\end{equation}
Then,
\begin{equation}
(R_I,R_J)=(-1)^{\ell(I)-1} \delta_{I,\bar J^\sim}
\end{equation}
so that this is (up to an unessential sign) the Bergeron-Zabrocki scalar
product~\cite[Eq.~(4)]{BZ}.
Indeed, if $\Des(I)=\{d_1,\dots,d_r\}$ and $\Des(J)=\{e_1,\dots,e_s\}$,
then
\begin{equation}
R_I^* R_J=(-1)^{d_r}\eta_{d_r}\dots (-1)^{d_1}\eta_{d_1}
\,\eta_{e_1}\dots\eta_{e_s}
%\hat{\eta}_{d_1}\dots\hat{\eta}_{d_r}\dots\eta_{e_s}
\end{equation}
and the coefficient of $\eta_1\dots\eta_{n-1}$ in this product is
zero if $\Des(I)$ and $\Des(J)$ are not complementary subsets of
$[n-1]$. When it is the case, moving $\eta_{d_k}$ to its place
in the middle of the $e_i$ produces a sign $(-1)^{d_k-1}$,
which together with the factor $(-1)^{d_k}$ results into a single
factor $(-1)$. Hence the final sign $(-1)^r=(-1)^{\ell(I)-1}$. 

%%%%%%%%%%%%%%%%%%%%%%%%%%%%%%%%%%%%%%%%%%%%%%%%%%%%%%%%
\subsection{Factorized elements in the Grassman algebra}

Now, for a sequence of parameters $Z=(z_1,\dots,z_{n-1})$, let
\begin{equation}
K_n(Z)=(1+z_1\eta_1)(1+z_2\eta_2) \dots (1+z_{n-1}\eta_{n-1})\,.
\end{equation}

Note that this is equivalent to define, as was already done
in~\cite[Eq.~(18)]{HLT},
\begin{equation}
K_n(A;Z) = \sum_{|I|=n} \left(\prod_{d\in \Des(I)} z_d\right) R_I(A) \,.
\end{equation}

For example, with $n=4$, if one orders compositions as usual by reverse
lexicographic order, the coefficients of expansion of $K_n$ on the  basis $(R_I)$
are
\begin{equation}
1, \ \ z_3, \ \ z_2, \ \ z_2z_3, \ \ z_1, \ \ z_1z_3, \ \ z_1z_2, \ \ z_1z_2z_3.
\end{equation}

We then have
\begin{lemma}
\begin{equation}
(K_n(X),K_n(Y))=\prod_{i=1}^{n-1}(y_i-x_i)\,.
\end{equation}
\end{lemma}

\Proof By induction. For $n=1$ the scalar product is $1$, 
and
\begin{eqnarray}
(K_{n+1}(X),K_{n+1}(Y)) &=& 
{\int} d\eta\,
(1+(-1)^{n}x_{n}\eta_n)(K_n(X),K_n(Y))\eta_1\dots\eta_{n-1}
(1+y_n\eta_n)\nonumber\\
&=&
{\int} d\eta\,
(K_n(X),K_n(Y))\eta_1\dots\eta_{n-1}(1-x_n\eta_n)(1+y_n\eta_n)\nonumber \\
&=&(y_n-x_n)(K_n(X),K_n(Y))\,.
\end{eqnarray}  
\qed

%%%%%%%%%%%%%%%%%%%%%%%%%%%%%%%%%%%%%%%%%%%%%%%%%%%%%%%%%%%%%%%%%%%%%%%%%%%%%%
\subsection{Bases of $\NCSF$}

We shall investigate bases of $\Sym_n$ of the form 
\begin{equation}
\label{bases}
\tHH_I=K_n(Z_I) = \sum_J \bkt_{IJ}R_J\,,
\end{equation}
where $Z_I$ is a sequence of parameters depending on the composition $I$ of
$n$.

The bases defined in~\cite{HLT} and~\cite{BZ} are of the previous form
and for both of them, the determinant of the Kostka matrix $\KK=(\bkt_{IJ})$
is a product of linear factors (as for ordinary Macdonald polynomials).
This is explained by the fact that these matrices have the form
\begin{equation}
\label{recmat}
\left(
\begin{matrix} 
 A & xA \\
 B & yB
\end{matrix}
\right)
\end{equation}
where $A$ and $B$ have a similar structure, and so on recursively.
Indeed, for such matrices,

\begin{lemma}
\label{lemdet}
Let $A,B$ be two $m\times m$ matrices. Then,
\begin{equation}
\left|
\begin{matrix} 
 A & xA \\
 B & yB
\end{matrix}
\right|
= (y-x)^m \det A\cdot \det B \,.
\end{equation}
\end{lemma}

%\Proof
%Subtract $x$ times the $j$th column to the $m+j$th column
%for all $1\le j\le m$.
%\hfill\qed

%%%%%%%%%%%%%%%%%%%%
\subsection{Duality}

Similarly, the dual vector space $\QSym_n=\Sym_n^*$ can be identified with a
Grassmann algebra on another set of generators $\xi_1,\dots,\xi_{n-1}$.
Encoding the fundamental basis $F_I$ of Gessel~\cite{Ges} by 
\begin{equation}
\xi_D:=\xi_{d1}\xi_{d_2}\dots \xi_{d_k},
\end{equation} 
the usual duality pairing such that the $F_I$ are dual to the $R_I$ is given
in this setting by
\begin{equation}
\<\xi_D,\eta_E\>=\delta_{DE}\,.
\end{equation}

Let
\begin{equation}
L_n(Z)=(z_1-\xi_1)\dots (z_{n-1}-\xi_{n-1})\,.
\end{equation}

Then, as above, we have a factorization identity:
\begin{lemma}
\begin{equation}
\<L_n(X),K_n(Y)\>=\prod_{i=1}^{n-1}(x_i-y_i)\,.
\end{equation}
\end{lemma}

\Proof By definition
\begin{equation}
L_n(X)=\sum_{D\subseteq[n-1]}(-1)^{|D|}\xi_D\prod_{e\not\in D}x_e
\end{equation}
and
\begin{equation}
K_n(Y)=\sum_{D\subseteq[n-1]}\eta_D\prod_{d\in D}y_d\,,
\end{equation}
so that
\begin{equation}
\<L_n(X),K_n(Y)\>
=\sum_{D\subseteq[n-1]}(-1)^{|D|}
\prod_{e\not\in D}x_e
\prod_{d\in D}y_d
=\prod_{i=1}^{n-1}(x_i-y_i)\,.
\end{equation}
\hfill\qed

Note that alternatively, assuming that the $\xi_i$ and the $\eta_j$ commute
with each other and that $\xi_i\eta_i=1$, one can define $\<f,g\>$ as the
constant term in the product $fg$.

Using this formalism, one can for example find for the dual basis
$\Phi_I^*$ of $\Phi^I$, the following expression
\begin{equation}
\Phi_J^* \longleftrightarrow
\prod_{i=1}^r \frac{1}{j_i}
\prod_{k\not\in \Des(K)}^\rightarrow (1-\xi_k)
     f_{r-1}(\xi_{d_1},\dots,\xi_{d_{r-1}}),
\end{equation}
where $\Phi^*_{1^r} \Longleftrightarrow f_{r-1}(\xi_1,\dots,\xi_{r-1})$,
which is simpler than the description of~\cite[Prop. 4.29]{NCSF1}.
Moreover, one can show that
\begin{equation}
\Phi_{1^n}^*(X) = F_n(X{\mathbb E})
                = \sum_{I\vdash n} R_I({\mathbb E}) F_I(X),
\end{equation}
where ${\mathbb E}$ is the exponential alphabet defined by
$S_n({\mathbb E}) = 1/n!$, so that
\begin{equation}
f_{r-1}(\xi_1,\ldots,\xi_{r-1}) = \sum_{D\subseteq [r-1]} a_D \xi_D,
\end{equation}
where $a_D$ is the number of permutations of $\SG_r$ with descent set $D$.

%%%%%%%%%%%%%%%%%%%%%%%%%%%%%%%%%%%%%%%%%%%%%%%%%%%%%%%%%%%%%%%%%%%%%%%%%%%%%%
%%%%%%%%%%%%%%%%%%%%%%%%%%%%%%%%%%%%%%%%%%%%%%%%%%%%%%%%%%%%%%%%%%%%%%%%%%%%%%
%%%%%%%%%%%%%%%%%%%%%%%%%%%%%%%%%%%%%%%%%%%%%%%%%%%%%%%%%%%%%%%%%%%%%%%%%%%%%%
\section{Bases associated with paths in a binary tree}

The most general possibility to build bases whose Kostka matrix can be
recursively decomposed into blocks of the form~(\ref{recmat}) is as follows.
Let $ \y= \{ y_u \}$ be a family of indeterminates indexed by all boolean
words of length $\leq n-1$. For example, for $n=3$, we have the six parameters
$y_0, y_1, y_{00}, y_{01}, y_{10}, y_{11}$.

We can encode a composition $I$ with descent set $D$ by the boolean word
$u=(u_1,\dots,u_{n-1})$ such that $u_i=1$ if $i\in D$ and $u_i=0$
otherwise.

Let us denote by $u_{m\dots p}$ the sequence $u_{m}u_{m+1}\dots u_{p}$
and define
\begin{equation}
 P_I := (1+ y_{u_1} \eta_1)
        (1+y_{u_{1\dots2}} \eta_2) \dots
        (1+y_{u_{1\dots n-1}} \eta_{n-1})
\end{equation}
or, equivalently,
\begin{equation}
\label{defP}
 P_I := K_n(Y_I)
\qquad\text{with}\quad Y_I=[ y_{u_1},y_{u_{1\dots2}},\dots,y_{u}] =:(y_k(I))\,.
\end{equation}

Similarly, let
\begin{equation}
Q_I := (y_{w_1} - \xi_1) (y_{w_{1\dots2}}- \xi_2) \dots
       (y_{w_{1\dots n-1}}- \xi_{n-1}) =: L_n(Y^I)
\end{equation}
where 
$ w_{1\dots k} = u_1\dots u_{k-1}\,\overline{u_k}$ where
$\overline{u_k}=1-u_k$, so that
\begin{equation}
Y^I=[ y_{w_1},y_{w_{1\dots2}},\dots,y_{w_{1\dots n-1}}]=(y^k(I))\,.
\end{equation}

For $n=4$, we have the following tables:
\small
\begin{equation}
\begin{array}{|c|c|c||c|c|}
\hline
 I  &  u  &   P_I 
&  Q_I \\
\hline
4   & 000 & (1+y_0\eta_1)(1+y_{00}\eta_2)(1+y_{000}\eta_3)
& (y_1-\xi_1)(y_{01}-\xi_2)(y_{001}-\xi_3)\\
31  & 001 & (1+y_0\eta_1)(1+y_{00}\eta_2)(1+y_{001}\eta_3)
& (y_1-\xi_1)(y_{01}-\xi_2)(y_{000}-\xi_3)\\
22  & 010 & (1+y_0\eta_1)(1+y_{01}\eta_2)(1+y_{010}\eta_3)
& (y_1-\xi_1)(y_{00}-\xi_2)(y_{011}-\xi_3)\\
211 & 011 & (1+y_0\eta_1)(1+y_{01}\eta_2)(1+y_{011}\eta_3)
& (y_1-\xi_1)(y_{00}-\xi_2)(y_{010}-\xi_3)\\
13  & 100 & (1+y_1\eta_1)(1+y_{10}\eta_2)(1+y_{100}\eta_3)
& (y_0-\xi_1)(y_{11}-\xi_2)(y_{101}-\xi_3)\\
121 & 101 & (1+y_1\eta_1)(1+y_{10}\eta_2)(1+y_{101}\eta_3)
& (y_0-\xi_1)(y_{11}-\xi_2)(y_{100}-\xi_3)\\
112 & 110 & (1+y_1\eta_1)(1+y_{11}\eta_2)(1+y_{110}\eta_3)
& (y_0-\xi_1)(y_{10}-\xi_2)(y_{111}-\xi_3)\\
1111& 111 & (1+y_1\eta_1)(1+y_{11}\eta_2)(1+y_{111}\eta_3)
& (y_0-\xi_1)(y_{10}-\xi_2)(y_{110}-\xi_3)\\
\hline
\end{array}
\end{equation}

\normalsize

%%%%%%%%%%%%%%%%%%%%%%%%%%%%%%%%%%%%%%%%%%%%%%%%%%%%%%%%%%%%%%%%%%%%%%%%%%%%%%
\subsection{Kostka matrices}

The Kostka matrix, which is defined as the transpose of the transition matrices from $P_I$ to $R_J$, 
is, for $n=4$:
\begin{equation}
\label{PtoR}
\left(
\begin{array}{cccccccc}
1 & y_{000} & y_{00} & y_{00}y_{000} & y_0 & y_0y_{000} & y_0y_{00}
  & y_0y_{00}y_{000} \\
1 & y_{001} & y_{00} & y_{00}y_{001} & y_0 & y_0y_{001} & y_0y_{00}
  & y_0y_{00}y_{001} \\
1 & y_{010} & y_{01} & y_{01}y_{010} & y_0 & y_0y_{010} & y_0y_{01}
  & y_0y_{01}y_{010} \\
1 & y_{011} & y_{01} & y_{01}y_{011} & y_0 & y_0y_{011} & y_0y_{01}
  & y_0y_{01}y_{011} \\
1 & y_{100} & y_{10} & y_{10}y_{100} & y_1 & y_1y_{100} & y_1y_{10}
  & y_1y_{10}y_{100} \\
1 & y_{101} & y_{10} & y_{10}y_{101} & y_1 & y_1y_{101} & y_1y_{10}
  & y_1y_{10}y_{101} \\
1 & y_{110} & y_{11} & y_{11}y_{110} & y_1 & y_1y_{110} & y_1y_{11}
  & y_1y_{11}y_{110} \\
1 & y_{111} & y_{11} & y_{11}y_{111} & y_1 & y_1y_{111} & y_1y_{11}
  & y_1y_{11}y_{111} \\
\end{array}
\right)
\end{equation}
For example,
\begin{equation}
\begin{split}
P_{211}
& = R_4 + y_{011}R_{31} + y_{01}R_{22} + y_{01}y_{011}R_{211} \\
& + y_0R_{13} + y_0y_{011}R_{121} + y_0y_{01}R_{112}
         + y_0y_{01}y_{011}R_{1111}.
\end{split}
\end{equation}
Note that this matrix is recursively of the form of Eq.~(\ref{recmat}).
Thus, its determinant is
\begin{equation}
(y_1-y_0)^4 (y_{01}-y_{00})^2 (y_{11}-y_{10})^2
(y_{001}-y_{000}) (y_{011}-y_{010}) (y_{101}-y_{100}) (y_{111}-y_{110}).
\end{equation}
This has the consequence that, given a specialization of the $y_w$, one can easily
check whether the $P_I$ form a linear basis of $\NCSF_n$.

\begin{proposition}
\label{dualPQ}
The bases $(P_I)$ and $(Q_I)$ are adjoint to each other, up to normalization:
\begin{equation}
\<Q_I,P_J\>=\<L_n(Y^I),K_n(Y_J)\>= \prod_{k=1}^{n-1}(y^k(I)-y_k(J)\>\,,
\end{equation}
which is indeed zero unless $I=J$.
\end{proposition}

\Proof If $I=J$, then $y^k(I)\not=y_k(I)$ by definition. If $I\not=J$,
let $d$ be the smallest integer which is a descent of either $I$
or $J$ but not both. Then, $y_{w_{1\dots d}}=y_{u_{1\dots d}}$
and $\<Q_I,P_J\>=0$.
\qed

From this, it is easy to derive a product
formula for the basis $P_I$. Note that we are considering the usual product
from $\Sym_n\times\Sym_m$ to $\Sym_{n+m}$ and not the product of the
Grassman algebra. 

\begin{proposition}
\label{prop-pipj}
Let $I$ and $J$ be two compositions of respective sizes $n$ and $m$.
The product $P_IP_J$ is a sum over an interval of the lattice of compositions
\begin{equation}
\label{pipj}
P_I P_J =\sum_{K\in [I\triangleright (m),I\cdot (1^m)]} c_{IJ}^K P_K
\end{equation}
where
\begin{equation}
\label{cijk}
c_{IJ}^K=\frac{\<L_{n+m}(Y^K),K_{n+m}(Y_I+1+Y_J)\>}{\<Q_K,P_K \>},
\end{equation}
where $Y_I+1+Y_J$ stands for the sequence
$(y_1(I),\dots,y_n(I),1,y_1(J),\dots,y_m(J))$.
\end{proposition}

\Proof
The usual product from $\Sym_n\times\Sym_m$ to $\Sym_{n+m}$ can be expressed
as
\begin{equation}
f_n\times g_m
 = f_n(\eta_1,\dots,\eta_{n-1})(1+\eta_n)g_m(\eta_{n+1},\dots,\eta_{n+m-1})\,.
\end{equation}
Indeed, this formula is clearly satisfied for $R_IR_J$. Then,
\begin{equation}
P_I P_J
 = K_n(\eta_1,\dots,\eta_{n-1})(1+\eta_n)K_m(\eta_{n+1},\dots,\eta_{n+m-1})
 = K_{n+m}(Y_I+1+Y_J).
\end{equation}
Then, for any $K$, the coefficient of $P_K$ is given by Formula~(\ref{cijk}).
Thanks to Proposition~\ref{dualPQ}, it is $0$
if the boolean vector corresponding to $I$ is not a prefix of the boolean
vector corresponding to $K$, that is, if $K$ is not in the interval
$[I\triangleright (m),I\cdot (1^m)]$.
\qed 

For example,
\begin{equation}
\label{p2p2}
\begin{split}
P_{2} P_{2}
& =
\frac{(y_{01}-1)(y_{001}-y_0)}{(y_{01}-y_{00})(y_{001}-y_{000})} P_{4} +
\frac{(y_{01}-1)(y_{000}-y_0)}{(y_{01}-y_{00})(y_{000}-y_{001})} P_{31} \\
& +
\frac{(y_{00}-1)(y_{011}-y_0)}{(y_{00}-y_{01})(y_{011}-y_{010})} P_{22} +
\frac{(y_{00}-1)(y_{010}-y_0)}{(y_{00}-y_{01})(y_{010}-y_{011})} P_{211}.
\end{split}
\end{equation}
\begin{equation}
\label{p1p1}
\begin{split}
P_{11} P_{11}
& =
\frac{(y_{11}-1)(y_{101}-y_1)}{(y_{11}-y_{10})(y_{101}-y_{100})} P_{13} +
\frac{(y_{11}-1)(y_{100}-y_1)}{(y_{11}-y_{10})(y_{100}-y_{101})} P_{121} \\
& +
\frac{(y_{10}-1)(y_{111}-y_1)}{(y_{10}-y_{11})(y_{111}-y_{110})} P_{112} +
\frac{(y_{10}-1)(y_{110}-y_1)}{(y_{10}-y_{11})(y_{110}-y_{111})} P_{1111}.
\end{split}
\end{equation}

%%%%%%%%%%%%%%%%%%%%%%%%%%%%%%%%%%%%%%%%%%%%%%%%%%%%%%%%%%%%%%%%%%%%%%%%%%%%%%
\subsection{The quasi-symmetric side}

As we have seen before, the $(Q_I)$ being dual to the $(P_I)$, the inverse
Kostka matrix is given by the simple construction:

%Define $Z'(I)=Z({}^tQ,{}^tT)(I)$ and $Z''(Q,T)=Z'(T,Q)$. 

\begin{proposition}
\label{invKgen}
The inverse of the Kostka matrix is given by
\begin{equation}
(\KK_n^{-1})_{IJ}
 =(-1)^{\ell(I)-1}
  \prod_{d\in\Des(\bar I^\sim)} y^d(J)
  \, \prod_{p=1}^{n-1} \frac1{y^p(J)-y_p(J)}\,.
\end{equation}
\end{proposition}

\Proof
This follows from Proposition~\ref{dualPQ}.
\qed

One can check the answer on the table for $n=4$ of $Q_I$ in terms of the
$F_I$.

%%%%%%%%%%%%%%%%%%%%%%%%%%%%%%%%%%%%%%%%%%%%%%%%%%%%%%%%%%%%%%%%%%%%%%%%%%%%%%
\subsection{Some specializations}

Let us now consider the specialization sending all $y_w$ to $1$ if $w$ ends
with a $1$ and denote by $\KK'$ the matrix obtained by this specialization.
Then, as in~\cite[p.~10]{HLT},

\begin{proposition}
\label{t1q1}
Let $n$ be an integer. Then
\begin{equation}
S_n = \KK_n {\KK'_n}^{-1}
\end{equation}
is lower triangular.
More precisely,
let $Y'_J$ be the image of $Y_J$ by the previous specialization and define
$Y'^J$ in the same way.
Then the coefficient $s_{IJ}$ indexed by $(I,J)$ is
\begin{equation}
s_{IJ} = \prod_{k=1}^{n-1} \frac{y_k(I)-y'^k(J)}{y'_k(J)-y'^k(J)}.
\end{equation}
\end{proposition}

\Proof
This follows from the explicit form of $\KK_n$ and ${\KK'_n}^{-1}$
(see Proposition~\ref{invKgen}).
\qed

Similarly, the specialization sending all $y_w$ to $1$ if $w$ ends with a $0$
leads to an upper triangular matrix.

These properties can be regarded as analogues of the characterization of
Macdonald polynomials given in~\cite[Prop.~2.6]{Hai}
(see also~\cite[p. 10]{HLT}).

%%%%%%%%%%%%%%%%%%%%%%%%%%%%%%%%%%%%%%%%%%%%%%%%%%%%%%%%%%%%%%%%%%%%%%%%%%%%%%
%%%%%%%%%%%%%%%%%%%%%%%%%%%%%%%%%%%%%%%%%%%%%%%%%%%%%%%%%%%%%%%%%%%%%%%%%%%%%%
%%%%%%%%%%%%%%%%%%%%%%%%%%%%%%%%%%%%%%%%%%%%%%%%%%%%%%%%%%%%%%%%%%%%%%%%%%%%%%
\section{The two-matrix family}

%%%%%%%%%%%%%%%%%%%%%%%%%%%%%%%%%%%%%%%%%%%%%%%%%%%%%%%%%%%%%%%%%%%%%%%%%%%%%%
\subsection{A specialization of the paths in a binary tree}

The above bases can now be specialized to bases $\tHH(A;Q,T)$, depending on
two infinite matrices of parameters.
Label the cells of the ribbon diagram of a composition $I$ of $n$ with their
matrix coordinates as follows:

\small
\begin{equation}
{\rm Diagr\,}(4,1,2,1)=
\moyyoung{(1,1) & (1,2)   & (1,3)   & (1,4) \cr
         \blank & \blank  & \blank  & (2,4) \cr
         \blank & \blank  & \blank  & (3,4) & (3,5)\cr
         \blank & \blank  & \blank  &\blank & (4,5)\cr}
\end{equation}
\normalsize
We associate a variable $z_{ij}$ with each cell except $(1,1)$ by setting
$z_{ij}:=q_{i,j-1}$ if $(i,j)$ has a cell on its left,
and $z_{ij}:=t_{i-1,j}$ if $(i,j)$ has a cell on its top.
The alphabet $Z(I)=(z_j(I))$ is the sequence of the $z_{ij}$ in their
natural order. For example,
\small
\begin{equation}
{Z}(4,1,2,1)=
\moyyoung{ & q_{11}   & q_{12}   & q_{13} \cr
         \blank & \blank  & \blank  & t_{14} \cr
         \blank & \blank  & \blank  & t_{24} & q_{34}\cr
         \blank & \blank  & \blank  &\blank & t_{35}\cr}
= ( q_{11},q_{12}, q_{13}, t_{14}, t_{24,} q_{34}, t_{35})
\end{equation}
\normalsize

Next, if $J$ is a composition of the same integer $n$, form the
monomial
\begin{equation}
\bkt_{IJ}(Q,T)=\prod_{d\in\Des(J)}z_d(I)\,.
\end{equation}
For example, with $I=(4,1,2,1)$ and $J=(2,1,1,2,2)$,
we have $\Des(J)=\{2,3,4,6\}$ and $\bkt_{IJ}=q_{12}q_{13}t_{14}q_{34}$.

\begin{definition}
Let $Q=(q_{ij})$ and $T=(t_{ij})$ ($i,j\ge 1$) be two infinite matrices
of commuting indeterminates.
For a composition $I$ of $n$, the noncommutative $(Q,T)$-Macdonald polynomial
$\tHH_I(A;Q,T)$ is
\begin{equation}
\label{defHt}
\tHH_I(A;Q,T)= K_n(A;Z(I)) = \sum_{J\vDash n}\bkt_{IJ}(Q,T)R_J(A)\,.
\end{equation}
Note that $\tHH_I$ depends only on the $q_{ij}$ and $t_{ij}$ with $i+j\le n$.
\end{definition}

\begin{note}
\label{rem}
{\rm
Since the $k$th element of $Z_I$ depends only on the prefix of size $k$ of the
boolean vector associated with $I$, the $\tHH_I$ are specializations of the
$P_I$ defined in Equation~(\ref{defP}).
More precisely, if $u=w0$ is a binary word ending by $0$, the specialization
is $y_u= q_{|w|_1+1,|w|_0+1}$, and if $u=w1$, we have
$y_u=t_{|w|_1+1,|w|_0+1}$.
Note also that $y_{w0}$ and $y_{w1}$ where $w$ is any binary word are
different (one is a $q$, the other is a $t$), so that the determinant of this
specialized Kostka matrix is generically nonzero.
Finally,  since these $\tHH$ are specializations of the $P$s,
the product formula detailed in Proposition~\ref{prop-pipj} gives a simple
generic product formula for the $\tHH$.
}
\end{note}

For example, translating Equations~(\ref{p2p2}) and~(\ref{p1p1}), one gets
\begin{equation}
\label{ht2ht2}
\begin{split}
\tHH_{2} \tHH_{2}
& =
\frac{(t_{12}-1)(t_{13}-q_{11})}{(t_{12}-q_{12})(t_{13}-q_{13})} \tHH_{4} +
\frac{(t_{12}-1)(q_{13}-q_{11})}{(t_{12}-q_{12})(q_{13}-t_{13})} \tHH_{31} \\
& +
\frac{(q_{12}-1)(t_{22}-q_{11})}{(q_{12}-t_{12})(t_{22}-q_{22})} \tHH_{22} +
\frac{(q_{12}-1)(q_{22}-q_{11})}{(q_{12}-t_{12})(q_{22}-t_{22})} \tHH_{211}.
\end{split}
\end{equation}
\begin{equation}
\label{ht1ht1}
\begin{split}
\tHH_{11} \tHH_{11}
& =
\frac{(t_{21}-1)(t_{22}-t_{11})}{(t_{21}-q_{21})(t_{22}-q_{22})} \tHH_{13} +
\frac{(t_{21}-1)(q_{22}-t_{11})}{(t_{21}-q_{21})(q_{22}-t_{22})} \tHH_{121} \\
& +
\frac{(q_{21}-1)(t_{31}-t_{11})}{(q_{21}-t_{21})(t_{31}-q_{31})} \tHH_{112} +
\frac{(q_{21}-1)(q_{31}-t_{11})}{(q_{21}-t_{21})(q_{31}-t_{31})} \tHH_{1111}.
\end{split}
\end{equation}

\subsection{$(Q,T)$-Kostka matrices}

Here are the $(Q,T)$-Kostka matrices for $n=3,4$ (compositions
are as usual in reverse lexicographic order):

\begin{equation}
\KK_3 =
\left(\begin{array}{cccc}
1 & q_{12} & q_{11} & q_{11} q_{12}\\
1 & t_{12} & q_{11} & q_{11} t_{12}\\
1 & q_{21\ } & t_{11} & q_{21} t_{11}\\
1 & t_{21} & t_{11} & t_{11} t_{21}
\end{array}\right)
\end{equation}

\begin{equation}
\KK_4=
\left(\begin{array}{cccccccc}
1 & q_{13} & q_{12} & q_{12} q_{13} & q_{11} & q_{11} q_{13} & q_{11} q_{12} & q_{11} q_{12} q_{13}\\
1 & t_{13} & q_{12} & q_{12} t_{13} & q_{11} & q_{11} t_{13} & q_{11} q_{12} & q_{11} q_{12} t_{13}\\
1 & q_{22} & t_{12} & q_{22} t_{12} & q_{11} & q_{11} q_{22} & q_{11} t_{12} & q_{11} q_{22} t_{12}\\
1 & t_{22} & t_{12} & t_{12} t_{22} & q_{11} & q_{11} t_{22} & q_{11} t_{12} & q_{11} t_{12} t_{22}\\
1 & q_{22} & q_{21} & q_{21} q_{22} & t_{11} & q_{22} t_{11} & q_{21} t_{11} & q_{21} q_{22} t_{11}\\
1 & t_{22} & q_{21} & q_{21} t_{22} & t_{11} & t_{11} t_{22} & q_{21} t_{11} & q_{21} t_{11} t_{22}\\
1 & q_{31} & t_{21} & q_{31} t_{21} & t_{11} & q_{31} t_{11} & t_{11} t_{21} & q_{31} t_{11} t_{21}\\
1 & t_{31} & t_{21} & t_{21} t_{31} & t_{11} & t_{11} t_{31} & t_{11} t_{21} & t_{11} t_{21} t_{31}
\end{array}\right)
\end{equation}

The factorization property of the determinant of the $(Q,T)$-Kostka 
matrix, which is valid for the usual Macdonald polynomials as well as
for the noncommutative analogues of \cite{HLT} and \cite{BZ} still holds since
the $\tHH_I$ are specializations of the $P_I$. More precisely,

\begin{theorem}
Let $n$ be an integer. Then
\begin{equation}
\det\KK_n = \prod_{i+j\le n}(q_{ij}-t_{ij})^{e(i,j)}\,,
\end{equation}
where $e(i,j)=\binom{i+j-2}{i-1}\, 2^{n-i-j}$.
\end{theorem}

\Proof
The matrix $\KK_n(Q,T)$ is of the form
\begin{equation}
\KK_n=\left(\begin{matrix}
A&q_{11}A\\ B&t_{11}B
\end{matrix}\right)
\end{equation}
where $A$ and $B$ are obtained from $\KK_{n-1}$ by the replacements
$q_{ij}\mapsto q_{i,j+1}$, $t_{ij}\mapsto t_{i,j+1}$,
and $q_{ij}\mapsto q_{i+1,j}$, $t_{ij}\mapsto t_{i+1,j}$, respectively.
So the result follows from Lemma~\ref{lemdet}.
\qed 

For example, with $n=4$, 
\begin{equation}
\det (\KK_4(Q,T))=
(q_{11} - t_{11})^4 (q_{12} - t_{12})^2  (q_{21} - t_{21})^2
(q_{22} - t_{22})^2  (q_{13} - t_{13})  (q_{31} - t_{31}).
\end{equation}

\subsection{Specializations}

For appropriate specializations, we recover (up to indexation) the Bergeron-Zabrocki
polynomials $\tHH_I^{BZ}$ of~\cite{BZ} and the multiparameter Macdonald
functions $\tHH_I^{HLT}$ of~\cite{HLT}:

\begin{proposition}
Let $(q_i)$, $(t_i)$, $i\ge 1$ be two sequences of indeterminates.
For a composition $I$ of $n$,

(i) Let $\nu$ be the anti-involution of
$\NCSF$ defined by $\nu(S_n)=S_n$.
Under the specialization $q_{ij}=q_{i+j-1}$, $t_{ij}=t_{n+1-i-j}$,
$\tHH_I(Q,T)$ becomes a multiparameter version of $i \nu(\tHH_I^{BZ})$, to which
it reduces under the further specialization $q_i=q^i$ and $t_i=t^i$.

(ii)  Under the specialization $q_{ij}=q_{j}$, $t_{ij}=t_{i}$,
$\tHH_I(Q,T)$ reduces to $\tHH_I^{HLT}$.
\end{proposition}

\Proof Equation~(\ref{defHt}) gives directly \cite[Eq.~(36)]{BZ}
under the specialization (i) and \cite[Eqs.~(2), (6)]{HLT} under
the specialization (ii). \qed

%%%%%%%%%%%%%%%%%%%%%%%%%%%%%%%%%%%%%%%%%%%%%%%%%%%%%%%%%%%%%%%%%%%%%%%%%%%%%%
\subsection{The quasi-symmetric side}

Families of $(Q,T)$-quasi-symmetric functions can now be defined by duality by
specialization of the $(Q_I)$ defined in the general case.
The dual basis of $(\tHH_J)$ in $QSym$ will be denoted by $(\tGG_I)$. We have
\begin{equation}
\tGG_I(X;Q,T)= \sum_J  \btg_{IJ}(q,t) F_J(X)
\end{equation}
where the coefficients are given by the transposed inverse
of the  Kostka matrix: $(\btg_{IJ})={}^t(\bkt_{IJ})^{-1}$.

%Define $Z'(I)=Z({}^tQ,{}^tT)(I)$ and $Z''(Q,T)=Z'(T,Q)$. 
Let $Z'(I)(Q,T)=Z(I)(T,Q)=Z(\bar I^\sim)(Q,T)$. 
Then, thanks to Proposition~\ref{invKgen} and to the fact that changing the
last bit of a binary word amounts to change a $q$ into a $t$, we have

\begin{proposition}
The inverse of the $(Q,T)$-Kostka matrix is given by
\begin{equation}
(\KK_n^{-1})_{IJ}
= (-1)^{\ell(I)-1}
  \prod_{d\in\Des(\bar I^\sim)} z'_d(J)
  \,\prod_{p=1}^{n-1}\frac1{z_p(J)-z'_p(J)}\,.
\end{equation}
\end{proposition}

Note that, as in the more general case of parameters indexed by binary words
(see Proposition~\ref{t1q1}),
if one specializes all $t$ (resp. all $q$) to $1$, one then gets lower (resp.
upper) triangular matrices with explicit coefficients, hence generalizing the
observation of~\cite{HLT}.

%%%%%%%%%%%%%%%%%%%%%%%%%%%%%%%%%%%%%%%%%%%%%%%%%%%%%%%%%%%%%%%%%%%%%%%%%%%%%%
%%%%%%%%%%%%%%%%%%%%%%%%%%%%%%%%%%%%%%%%%%%%%%%%%%%%%%%%%%%%%%%%%%%%%%%%%%%%%%
%%%%%%%%%%%%%%%%%%%%%%%%%%%%%%%%%%%%%%%%%%%%%%%%%%%%%%%%%%%%%%%%%%%%%%%%%%%%%%
\section{Multivariate BZ polynomials}

In this section, we restrict our attention to the multiparameter
version of the Bergeron-Zabrocki polynomials, obtained
by setting $q_{ij}=q_{i+j-1}$ and  $t_{ij}=t_{n+1-i-j}$ in degree~$n$.

%%%%%%%%%%%%%%%%%%%%%%%%%%%%%%%%%%%%%%%%%%%%%%%%%%%%%%%%%%%%%%%%%%%%%%%%%%%%%%
\subsection{Multivariate BZ polynomials}

As in the case of the two matrices of parameters, $Q$ and $T$, one can deduce
the product in the $\tHH$ basis by some sort of specialization of the general
case. However, since $t_{ij}$ specializes to another $t$ where $n$ appears,
one has to be a little more cautious to get the correct answer.

\begin{theorem}
\label{thm-bz}
Let $I$ and $J$ be two compositions of respective sizes $p$ and $r$.
Let us denote by $K=I.\bar J^\sim$ and $n=|K|=p+r$.
Then
\begin{equation}
\tHH_I \tHH_J =
\frac{(-1)^{\ell(I)+|J|}}{\prod_{k\in\Des(K)} (q_k - t_{n-k})}
\sum_{K'} 
  \prod_{k\in\Des(K)} (-1)^{\ell(K)} (z_{k}(K')-z'_k(K')),
 \tHH_{K'}
\end{equation}
where the sum is computed as follows.
Let $I'$ and $J'$ be the compositions such that $|I'|=|I|$ and either
$K'=I'\cdot J'$, or $K'=I'\triangleright J'$.
If $I'$ is not coarser than $I$ or if $J'$ is not finer than $J$, then
$\tHH(K')$ does have coefficient $0$.
Otherwise,
%$z_k(K')= Z^{BZ}(K)_k$ if $k$ is a descent of $K'$, or
%$z_k(K')= Z^{BZ}(\bar K^\sim)_k$ if not.
%Moreover, $z'_k(K')=(Z^{BZ}(I').1.Z^{BZ}(J'))_k$ if $k$ is a descent of $K'$, 
%and
%$z'_k(K')=(Z^{BZ}(\bar{I'}^\sim).1.Z^{BZ}(\bar{J'}^\sim))_k$
%if not.
$z_k(K')=q_{k}$ if $k$ is a descent of $K'$ and $t_{n-k}$ otherwise.
Finally, $z'_k(K')$ does not depend on $K'$ and is
$(Z(I),1,Z(J))$.
\end{theorem}

For example,
with $I=J=(2)$, we have $K=(211)$.
Note that $Z^{BZ}(K)=[q_1,t_2,t_1]$ and
$Z^{BZ}(\bar K^\sim)=[t_3,q_2,q_3]$.
The set of compositions $K'$ having a nonzero coefficient is $(4)$, $(31)$,
$(22)$, $(211)$.
Here are the (modified) $Z$ and $Z'$ restricted to the descents of $K$
for these four compositions.
\begin{equation}
\begin{array}{lllll}
   & (4)        & (31)      & (22)      & (211)     \\
Z  & [t_2,t_1]  & [t_2,q_3] & [q_2,t_1] & [q_2,q_3] \\
Z' & [1,q_1]    & [1,q_1]   & [1,q_1]   & [1,q_1]   \\
\end{array}
\end{equation}

\begin{equation}
\label{bz22}
\begin{split}
\tHH^{BZ}_{2} \tHH^{BZ}_{2}
& =
\frac{(t_2-1)(t_1-q_1)}{(t_2-q_2)(t_1-q_3)} \tHH^{BZ}_{4} +
\frac{(t_2-1)(q_3-q_1)}{(t_2-q_2)(q_3-t_1)} \tHH^{BZ}_{31} \\
& +
\frac{(q_2-1)(t_1-q_1)}{(q_2-t_2)(t_1-q_3)} \tHH^{BZ}_{22} +
\frac{(q_2-1)(q_3-q_1)}{(q_2-t_2)(q_3-t_1)} \tHH^{BZ}_{211}.
\end{split}
\end{equation}

\begin{equation}
\label{bz1111}
\begin{split}
\tHH^{BZ}_{11} \tHH^{BZ}_{11}
& =
\frac{(t_2-1)(t_3-t_1)}{(t_2-q_2)(t_3-q_1)} \tHH^{BZ}_{31} +
\frac{(q_2-1)(t_3-t_1)}{(q_2-t_2)(t_3-q_1)} \tHH^{BZ}_{211} \\
& +
\frac{(t_2-1)(q_1-t_1)}{(t_2-q_2)(q_1-t_3)} \tHH^{BZ}_{121} +
\frac{(q_2-1)(q_1-t_1)}{(q_2-t_2)(q_1-t_3)} \tHH^{BZ}_{1111}.
\end{split}
\end{equation}

%%%%%%%%%%%%%%%%%%%%%%%%%%%%%%%%%%%%%%%%%%%%%%%%%%%%%%%%%%%%%%%%%%%%%%%%%%%%%%
\subsection{The $\nabla$ operator}

The $\nabla$ operator of~\cite{BZ} can be extended by
\begin{equation}
\label{eqnabla}
\nabla\tHH_I = \left(\prod_{d=1}^{n-1} z_d(I)\right) \,\,\tHH_I\,.
\end{equation}
Then,
\begin{proposition}
The action of $\nabla$ on the ribbon basis is given by
\begin{equation}
\nabla R_I=(-1)^{|I|+\ell(I)}
\prod_{d\in \Des(I^\sim)}q_d
\prod_{d\in \Des(\bar I^\sim)}t_d
\sum_{J\ge \bar I^\sim}\prod_{i\in\Des(I)\cap\Des(J)}(t_i+q_{n-i}) R_J\,.
\end{equation}
\end{proposition}

\Proof
This is a direct adaptation of the proof of \cite{BZ}.
Lemma 21, Corollary 22 and Lemma 23 of~\cite{BZ} remain valid
if one interprets $q^i$ and $t^j$ as $q_i$ and $t_j$.
In particular, Equation~(\ref{eqnabla}) reduces to \cite[(54)]{BZ}
under the specialization $q_i=q^i$, $t_i=t^i$.
\qed

Note also that if $I=(1^n)$, one has
\begin{equation}
\nabla \Lambda_n=\sum_{J\vDash n}\prod_{j\in\Des(J)}(q_j+t_{n-j}) R_J
= \sum_{J\vDash n}\prod_{j\not\in\Des(J)}(q_j+t_{n-j}-1)\Lambda^J\,.
\end{equation}

As a positive sum of ribbons, this is the multigraded characteristic 
of a projective module of the 0-Hecke algebra. Its dimension is the number of
packed words of length $n$ (called preference functions in \cite{BZ}).
Let us recall that a packed word is a word $w$ over $\{1,2,\dots\}$ so that if
$i>1$ appears in $w$, then $i-1$ also appears in $w$. The set of all packed
words of size $n$ is denoted by $\PW_n$.

Then the multigraded dimension of the previous module is
\begin{equation}
W_n(\q,\t)=\<\nabla\Lambda_n, F_1^n\>=\sum_{w\in\PW_n}\phi(w)
\end{equation}
where the statistic $\phi(w)$ is obtained as follows.

Let $\sigma_w=\overline{\std(\overline{w})}$, where $\overline{w}$ denotes
the mirror image of $w$.
Then
\begin{equation}
\phi(w)=\prod_{i\in\Des(\sigma_w^{-1})}x_i
\end{equation}
where $x_i=q_i$ if $w^\uparrow_i=w^\uparrow_{i+1}$ and $x_i=t_{n-i}$
otherwise, where $w^\uparrow$ is the nondecreasing reordering of $w$.

For example, with $w=22135411$, $\sigma_w=54368721$, $w^\uparrow=11122345$,
the recoils of $\sigma_w$ are $1$, $2$, $3$, $4$, $7$,
and $\phi(w)=q_1q_2t_5q_4t_1$.

Actually, we have the following slightly stronger result.

\begin{proposition}
Denote by $d_I$ the number of permutations $\sigma$ with descent composition
$C(\sigma)=I$. Then,
\begin{equation}
\label{nablam}
\nabla\Lambda_n
= \sum_{w\in\PW_n}
  \frac{\phi(w)}{d_{C(\sigma_w^{-1})}} R_{C(\sigma_w^{-1})}.
\end{equation}
\end{proposition}

\Proof
If $\sigma$ is any permutation such that $C(\sigma^{-1})=I$, the coefficient
of $R_I$ in the r.h.s. of~(\ref{nablam}) can be rewritten as
\begin{equation}
\label{nablam2}
\sum_{w\in\PW_n}
\frac{\phi(w)}{d_{C(\sigma_w^{-1})}}
=\sum_{\sigma_w=\sigma}\phi(w)\,.
\end{equation}

The words $w\in\PW_n$ such that $\sigma_w=\sigma$ are obtained by the
following construction.
For $i=\sigma^{-1}(1)$, we have $w_i=1$.
Next, if $k+1$ is to the right of $k$ in $\sigma$, we must have
$w_{\sigma^{-1}(k+1)}=w_{\sigma^{-1}(k)}+1$.
Otherwise, we have two choices:
$w_{\sigma^{-1}(k+1)}=w_{\sigma^{-1}(k)}$ or
$w_{\sigma^{-1}(k+1)}=w_{\sigma^{-1}(k)}+1$.
These choices are independent, and so contribute a factor $(q_i+t_{n-1})$ to
the sum~(\ref{nablam2}).
\qed

This can again be generalized:
\begin{theorem}
For any composition $I$ of $n$,
\begin{equation}\label{nablaRI}
\nabla R_I=(-1)^{|I|+\ell(I)} \theta(\sigma)
\sum_{w\in\PW_n;\,\ev(w)\le I}
\frac{R_{C(\sigma_w^{-1})}}{d_{C(\sigma_w^{-1})}}\,,
\end{equation}
where $\sigma$ is any permutation such that $C(\sigma^{-1})=\bar I^\sim$, and
\begin{equation}
\theta(\sigma)=\prod_{d\in\Des(\bar{I}^\sim)}t_d\,.
\end{equation}
\end{theorem}

\Proof
First, if $\ev(w)\le I$, then $C(\sigma_w^{-1})\ge \bar{I}^\sim$.
For any $J\ge \bar{I}^\sim$, the coefficient of $R_J$ in (\ref{nablaRI})
is, for any permutation $\tau$ such that $C(\tau^{-1})=J$,
\begin{equation}
\sum_{\sigma_w=\tau;\ \ev(w)\le I}\phi(w)\,.
\end{equation}
For a packed word $w$ such that $\sigma_w=\tau$, we have
\begin{equation}
\phi(w)=\prod_{j\in\Des(J)}x_j
=\prod_{j\in\Des(J)\cap\Des(I)}x_j
\prod_{j\in\Des(J)\backslash \Des(I)}x_j\,.
\end{equation}
In the second product, one has always $x_j=q_j$, since
$w^\uparrow_j=w^\uparrow_{j+1}$. In the first one,
there are, as before, two possible independent choices for each $j$.
\qed

The behaviour or the multiparameter BZ polynomials with respect to the scalar
product
\begin{equation}
[R_I,R_J]:= (-1)^{|I|+\ell(I)}\delta_{I,\bar J^\sim}
\end{equation}
is the natural generalization of \cite[Prop. 1.7]{BZ}:
\begin{equation}
[\tHH_I,\tHH_J]= (-1)^{|I|+\ell(I)}\delta_{I,\bar
J^\sim}\prod_{i=1}^{n-1}(q_i-t_{n-i})\,.
\end{equation}

%%%%%%%%%%%%%%%%%%%%%%%%%%%%%%%%%%%%%%%%%%%%%%%%%%%%%%%%%%%%%%%%%%%%%%%%%%%%%%%
%%%%%%%%%%%%%%%%%%%%%%%%%%%%%%%%%%%%%%%%%%%%%%%%%%%%%%%%%%%%%%%%%%%%%%%%%%%%%%%
%%%%%%%%%%%%%%%%%%%%%%%%%%%%%%%%%%%%%%%%%%%%%%%%%%%%%%%%%%%%%%%%%%%%%%%%%%%%%%%
\section{Quasideterminantal bases}

%%%%%%%%%%%%%%%%%%%%%%%%%%%%%%%%%%%%%%%%%%%%%%%%%%%%%%%%%%%%%%%%%%%%%%%%%%%%%%%
\subsection{Quasideterminants of almost triangular matrices}

Quasideterminants~\cite{GR} are noncommutative analogs of the ratio of a
determinant by one of its principal minors. Thus, the quasideterminants of a
generic matrix are not polynomials, but complicated rational expressions
living in the free field generated by the coefficients. However, for an almost
triangular matrices, \emph{i.e.}, such that $a_{ij}=0$ for $i>j+1$, all
quasideterminants are polynomials, with a simple explicit expression. We shall
only need the formula (see~\cite{NCSF1}, Prop.2.6):
\begin{equation}
\label{devel-qd}
\left|
\begin{matrix}
a_{11} & a_{12} & a_{13} & \dots & \bo{a_{1n}} \\
  -1   & a_{22} & a_{23} & \dots & a_{2n} \\
  0    &  -1    & a_{33} & \ddots& \vdots \\
 \vdots& \ddots & \ddots & \ddots& a_{n-1n} \\
   0   & \dots  & 0      &  -1   & a_{nn}
\end{matrix}
\right|
= a_{1n} +
 \sum_{1\leq j_1<\dots<j_k<n}
    a_{1j_1}a_{j_1+1j_2}a_{j_2+1 j_3} \dots a_{j_k+1 n}.
\end{equation}

Recall that the quasideterminant $|A|_{pq}$ is invariant by scaling the
rows of index different from $p$ and the columns of index diffrerent from $q$.
It is homogeneous of degree 1 with respect to row $p$ and column $q$. Also, the
quasideterminant is invariant under permutations of rows and columns.

The quasideterminant~(\ref{devel-qd}) coincides with the row-ordered expansion
of an ordinary determinant
\begin{equation}
\label{rowdet}
\rdet(A) := 
\sum_{\sigma\in{\SG_n}}
  \varepsilon(\sigma)
  a_{1\sigma(1)} a_{2\sigma(2)}\cdots a_{n\sigma(n)}
\end{equation}
which will be denoted as an ordinary determinant in the sequel.

%%%%%%%%%%%%%%%%%%%%%%%%%%%%%%%%%%%%%%%%%%%%%%%%%%%%%%%%%%%%%%%%%%%%%%%%%%%%%%%
\subsection{Quasideterminantal bases of $\Sym$}

Many interesting families of noncommutative symmetric
functions can be expressed as quasi-determinants of the form
\begin{equation}
H(W,G)=
\left|
\begin{matrix}
w_{11}G_1 & w_{12}G_2 & \dots & w_{1\, n-1}G_{n-1} & \bo{w_{1n}G_n} \\
w_{21}    & w_{22}G_1 & \dots & w_{2\, n-1}G_{n-2} & w_{2n}G_{n-1}  \\
0         & w_{32}    & \dots & w_{3\, n-3}G_{n-3} & w_{3n}G_{n-2}  \\
\vdots    & \vdots    & \ddots & \vdots             & \vdots         \\
0         & 0         & \dots & w_{n\, n-1}        & w_{nn}G_1      \\
\end{matrix}
\right|
\end{equation}
(or of the transposed form), where $G_k$ is some sequence of free generators
of $\Sym$, and $W$ an almost-triangular ($w_{ij}=0$ for $i>j+1$) scalar
matrix.
For example (see \cite[(37)-(41)]{NCSF1}),

\begin{equation}
S_n = (-1)^{n-1}
\left|\begin{matrix}
\Lambda_1 & \Lambda_2 & \dots & \Lambda_{n-1} &\bo{\Lambda_n}\\
\Lambda_0 & \Lambda_1 & \dots & \Lambda_{n-2} &\Lambda_{n-1} \\
0         & \Lambda_0 & \dots & \Lambda_{n-3} &\Lambda_{n-2} \\
\vdots    & \vdots    & \ddots & \vdots        &\vdots        \\
0         & 0         & \dots & \Lambda_0     &\Lambda_1     \\
\end{matrix}\right|\ ,
\end{equation}

\begin{equation} 
n\, S_n =
\left|\begin{matrix}
\Psi_1       & \Psi_2       & \dots & \Psi_{n-1} &\bo{\Psi_n}\\
-1           & \Psi_1       & \dots & \Psi_{n-2} &\Psi_{n-1} \\
0            & -2           & \dots & \Psi_{n-3} &\Psi_{n-2} \\
\vdots       & \vdots       & \ddots & \vdots     &\vdots     \\
0            & 0            & \dots & -n+1       &\Psi_1     \\
\end{matrix}\right|\ ,
\end{equation}
or (see~\cite[Eq.~(78)]{NCSF2})
\begin{equation}
[n]_q \, S_n(A)
=
\left|
\begin{matrix}
\Theta_1(q) & \Theta_2(q)    & \dots  & \Theta_{n-1}(q)
& \bo{\Theta_n(q)}  \\
-[1]_q      & q\,\Theta_1(q) & \dots  & q\,\Theta_{n-2}(q)
& q\,\Theta_{n-1}(q)\\ 
0           & - [2]_q        & \dots  & q^2\, \Theta_{n-3}(q)
& q^2\, \Theta_{n-2}(q) \\
\vdots      & \vdots         & \vdots & \vdots   & \vdots \\
0           & 0              & \dots  & -[n-1]_q & q^{n-1} \, \Theta_1(q)
\end{matrix}
\right|
\ ,
\end{equation}
where $\Theta_n(q)=(1-q)^{-1}S_n((1-q)A)$.
These examples illustrate relations between sequences of free generators.
Quasi-determinantal expressions for some linear bases can be recast
in this form as well. For example, the formula for ribbons
\begin{equation}
  (-1)^{n-1}R_I =
\left|
\begin{matrix}
S_{i_1} & S_{i_1+i_2} & S_{i_1+i_2+i_3} & \dots  & \bo{S_{i_1+\dots+i_n}} \\
S_0     & S_{i_2}     & S_{i_2+i_3}     & \dots  & S_{i_2+\dots+i_n} \\
0       & S_0         & S_{i_3}         & \dots  & S_{i_3+\dots+i_n} \\
\vdots  & \vdots      & \vdots          & \ddots & \vdots \\
0       & 0           & 0               & \dots  & S_{i_n}
\end{matrix}
\right|
\end{equation}
can be rewritten as follows. Let $U$ and $V$ be the $n\times n$
almost-triangular matrices

\begin{equation}
U=\left[\begin{matrix}
1 & 1 & \dots & 1 & 1 \\
-1 & -1 & \dots & -1 & -1 \\
0  & -1 & \dots & -1 & -1 \\
\vdots & \vdots & \ddots &\vdots &\vdots\\
0 & 0 & \dots & -1 & -1
\end{matrix}\right]
\quad
V=\left[\begin{matrix}
1 & 1 & \dots & 1 & 1 \\
-1 & 0 & \dots & 0 & 0 \\
0 & -1 & \dots & 0 & 0 \\
\vdots & \vdots & \ddots &\vdots &\vdots\\
0 & 0 & \dots & -1 & 0
\end{matrix}\right]
\end{equation}

Given the pair $(U,V)$, define, for each composition $I$ of $n$, a matrix
$W(I)$ by
\begin{equation}
w_{ij}(I)
=\begin{cases}
  u_{ij} & \text{if $i-1\in\Des(I)$},\\
  v_{ij} & \text{otherwise},
\end{cases}
\end{equation}
and set
\begin{equation}
H_I(U,V; A) := H(W(I),S(A))\,.
\end{equation}
Then,
\begin{equation}
(-1)^{\ell(I)-1} R_I = H_I(U,V)\,.
\end{equation}
Indeed, $H_I(U,V;A)$ is obtained by substituting in~(\ref{devel-qd})
\begin{equation}
a_{j_p+1,j_{p+1}}=
\begin{cases}
-S_{j_{p+1}-j_p} & \text{if $j_p\in\Des(I)$}\\
0 & \text{otherwise.}
\end{cases}
\end{equation}
This yields
\begin{equation}
\begin{split}
& S_n + \sum_k \sum_{\{j_1<\dots <j_k\}\subseteq
\Des(I)}S_{j_1}(-S_{j_2-j_1})\dots(-S_{n-j_k}) \\
&
=\sum_{\Des(K)\subseteq \Des(I)}(-1)^{\ell(K)-1}S^K=(-1)^{\ell(I)-1}R_I\,.
\end{split}
\end{equation}

For example,
\begin{equation}
R_{211}=
\left|\begin{matrix}
S_1 & S_2 &  S_3 & \bo{S_4} \\
-1  &  0  &   0  & 0        \\
0   &  -1 & -S_1 & -S_2     \\
0   &  0  &  -1  & -S_1
\end{matrix}\right|
= S_4 - S_{31} - S_{22} + S_{211}.
\end{equation}

For a generic pair of almost-triangular matrices $(U,V)$, the
$H_I$ form a basis of $\Sym_n$. Without loss of generality,
we may assume that $u_{1j}=v_{1j}=1$ for all $j$.
Then, the transition matrix $M$ expressing the $H_I$ on the
$S^J$ where $J=(j_1,\dots,j_p)$ satisfies:
\begin{equation}
M_{J,I} := x_{1 j_1-1} x_{j_1 j_2-1} \dots x_{j_p n}.
\end{equation}
where $x_{ij}=u_{ij}$ if $i-1$ is not a descent of $I$ and $v_{ij}$ otherwise.

As we shall sometimes need different
normalizations, we aslo define for arbitrary almost triangular matrices $U,V$
\begin{equation}\label{varH} % Attention au signe, choisir la bonne option ...
H'(W,G)= \rdet
\left[
\begin{matrix}
w_{11}G_1 & w_{12}G_2 & \ldots & w_{1\, n-1}G_{n-1} & {w_{1n}G_n} \\
w_{21}    & w_{22}G_1 & \ldots & w_{2\, n-1}G_{n-2} & w_{2n}G_{n-1}  \\
0         & w_{32}    & \ldots & w_{3\, n-3}G_{n-3} & w_{3n}G_{n-2}  \\
\vdots    & \vdots    & \ddots & \vdots             & \vdots         \\
0         & 0         & \ldots & w_{n\, n-1}        & w_{nn}G_1      \\
\end{matrix}
\right]
\end{equation}
and
\begin{equation}
H'_I(U,V) = H'(W(I),S(A))\,.
\end{equation}

%%%%%%%%%%%%%%%%%%%%%%%%%%%%%%%%%%%%%%%%%%%%%%%%%%%%%%%%%%%%%%%%%%%%%%%%%%%%%%%
\subsection{Expansion on the basis $(S^I)$}

For a composition $I=(i_1,\ldots,i_r)$ of $n$, let $I^\sharp$
be the integer vector of length $n$ obtained by replacing
each entry $k$ of $I$ by the sequence $(k,0,\ldots,0)$ ($k-1$
zeros):
\begin{equation}
I^\sharp = (i_10^{i_1-1}i_20^{i_2-1}\ldots i_r0^{i_r-1})\,,
\end{equation}
\emph{e.g.}, for compositions of 3,
\begin{equation}
(3)^\sharp=(300),\ (21)^\sharp=(201),\ (12)^\sharp=(120),\
(111)^{\sharp}=(111)\,.
\end{equation}
Adding (componentwise) the sequence $(0,1,2,\ldots,n-1)$ to $I^\sharp$,
we obtain a permutation $\sigma_I$. For example,
\begin{equation}
\sigma_{(3)}=(312),\ \sigma_{(21)}=(213),\ \sigma_{(12)}=(132),\
\sigma_{(111)}=(123)\,.
\end{equation}

\begin{proposition}
The expansion of $H'(W,S)$ on the $S$-basis is given by
\begin{equation}
H'(W,S)= \sum_{I\vDash n}\varepsilon(\sigma_I) w_{1\sigma_I(1)}\cdots
w_{n\sigma_I(n)} S^I. % signe global ?
\end{equation}
\end{proposition}

Thus, for $n=3$,
\begin{equation}
\begin{split}
H'(W,S)
=&
\left|
\begin{matrix}
w_{11} S_1 & w_{12} S_2 & w_{13} S_3 \\
w_{21}     & w_{22} S_1 & w_{23} S_2 \\
 0         & w_{32}     & w_{33} S_1 \\
\end{matrix}
\right| \\
=&
w_{11}w_{22}w_{33}S^{111}
-  w_{11}w_{23}w_{32}S^{12}
- w_{12}w_{21}w_{33}S^{21}
+ w_{13}w_{21}w_{32}S^{3}.
\end{split}
\end{equation}

%%%%%%%%%%%%%%%%%%%%%%%%%%%%%%%%%%%%%%%%%%%%%%%%%%%%%%%%%%%%%%%%%%%%%%%%%%%%%%
\subsection{Expansion of the basis $(R_I)$}

\begin{proposition}
For $I=(i_1,\ldots,i_r)$ be a composition of $n$, denote by $W_I$ the
product of diagonal minors of the matrix $W$ taken over the first $i_1$
rows and columns, then the next $i_2$ ones and so on. Then,
\begin{equation}
H'(W,S) = \sum_{I\vDash n}W_I R_I \,. % Gare aux signes !!
\end{equation}
\end{proposition}

%%%%%%%%%%%%%%%%%%%%%%%%%%%%%%%%%%%%%%%%%%%%%%%%%%%%%%%%%%%%%%%%%%%%%%%%%%%%%%%
\subsection{Examples}

\subsubsection{A family with factoring coefficients}

\begin{theorem}
Let $U$ and $V$ be defined by
\begin{equation}
u_{ij}=\begin{cases}
x^j-y^j                     & \text{if $i=1$}\\
aq_{i-1}x^{j-i+1}-y^{j-i+1} & \text{if $1<i<j+2$}\\
0                           & \text{otherwise}
\end{cases}
\end{equation}
\begin{equation}
%\quad
v_{ij}=\begin{cases}
x^j-y^j                        & \text{if $i=1$}\\
x^{j-i+1}-bu_{n+1-i}y^{j-i+1}  & \text{if $1<i<j+2$}\\
0                              & \text{otherwise}
\end{cases}
\end{equation}
Then the coefficients $W_J$ of the expansion of $H'_I(U,V)$
on the ribbon basis all factor as products of binomials.
\end{theorem}

\Proof Observe first that the substitutions $b=a^{-1}$ and
$u_{n+1-i}=q_{i-1}^{-1}$ changes the determinant of $V$ into
$c_n\det(U)$ with $c_n=a^{n-1}q_1\cdots q_{n-1}$. Expanding
$\det(U)$ by its first column yields a two-term recurrence 
implying easily the factorized expression
\begin{equation}
\det(U)=(x-y)\prod_{i=1}^{n-1}(x-aq_iy)
\end{equation}
so that as well
\begin{equation}
\det(V)=(x-y)\prod_{i=1}^{n-1}(y-bu_{n-i}x)\,.
\end{equation}
Now, all the matrices $W(I)$ built from $U$ and $V$, and all their diagonal
minors have the same structure, and their determinants factor similarly.
\qed

The formula for the coefficient of $R_n$ is simple enough: if one
orders the factors of $\det(U)$ and $\det(V)$ as
\begin{equation}
Z_n = (x-aq_1y, x-aq_2y,\ldots,x-aq_{n-1}y)
\end{equation}
and
\begin{equation}
Z'_n=(y-bu_{n-1}x,y-bu_{n-2}x,\ldots,y-bu_1x)
\end{equation}
then, the coefficient of $R_n$ in $H'_I(U,V)$
is
\begin{equation}
(x-y)\prod_{d\in\Des(I)}z'_d\prod_{e\not\in\Des(I)}z_e \,.
\end{equation} 

For example,

\begin{equation}
\begin{split}
\frac{H'_3(U,V)}{(x-y)}
 =& (x-aq_1y)(x-aq_2y) R_3
 + (x-aq_1y)(aq_2x-y) R_{21}\\
 &+ a(x-y)(q_1x-q_2y)  R_{12}
 + (aq_1x-y)(aq_2x-y) R_{111}.
\end{split}
\end{equation}

\begin{equation}
\begin{split}
\frac{H'_{21}(U,V)}{(x-y)}
 =& (x-aq_1y)(bu_1x-y) R_3
 + (x-aq_1y)(x-bu_1y) R_{21}\\
 &+ (abq_1u_1x-y)(x-y) R_{12}
 + (aq_1x-y)(x-bu_1y) R_{111}.
\end{split}
\end{equation}
\begin{equation}
\begin{split}
\frac{H'_{12}(U,V)}{(x-y)}
=& (x-aq_2y)(bu_2x-y) R_3
 + (aq_2x-y)(bu_2x-y) R_{21}\\
& + (x-abq_2u_2y)(x-y) R_{12}
 + (aq_2x-y)(x-bu_2y) R_{111}.
\end{split}
\end{equation}
\begin{equation}
\begin{split}
\frac{H'_{111}(U,V)}{(x-y)}
 =& (bu_1x-y)(bu_2x-y) R_3
 + (x-bu_1y)(bu_2x-y) R_{21}\\
 &+ b(u_1x-u_2y)(x-y)  R_{12}
 + (x-bu_1y)(x-bu_2y) R_{111}.
\end{split}
\end{equation}

\medskip

A more careful analysis allows one to
compute directly the 
coefficient of $R_J$ is $H'_I$:
denote by $u$ the boolean vector of $I$ and by $v$ the boolean vector of $J$
and consider the biword $w=\binom{u}{v}$.
Start with $c_{IJ}=c:=1$.
First, There are factors coming from the boundaries of the biword:
\begin{itemize}
\item If $w_1=\binom{0}{0}$ then $c:=(x-q_1y) c$, \\[-8pt]
\item If $w_1=\binom{1}{0}$ then $c:=(x-u_{n_1}y) c$, \\[-8pt]
\item If $w_{n-1}=\binom{0}{1}$ then $c:=(x-q_{n-1}y) c$, \\[-8pt]
\item If $w_{n-1}=\binom{1}{1}$ then $c:=(x-u_1y) c$.
\end{itemize}
Then, for any $i\in[1,n-2]$, the two biletters $w_i w_{i+1}$ can have
different values:
\begin{itemize}
\item If $w_i w_{i+1}=\binom{.\ 0}{0\ 0}$ then $c:=(x-q_{i+1}y) c$,\\[-8pt]
\item If $w_i w_{i+1}=\binom{.\ 1}{0\ 0}$ then $c:=(xu_{n-i-1}-y) c$,\\[-8pt]
\item If $w_i w_{i+1}=\binom{0\ .}{1\ 1}$ then $c:=(xq_i-y) c$,\\[-8pt]
\item If $w_i w_{i+1}=\binom{1\ .}{1\ 1}$ then $c:=(x-u_{n-i}y) c$,\\[-8pt]
\item If $w_i w_{i+1}=\binom{0\ 0}{1\ 0}$
      then $c:=(xq_i-q_{i+1}y)(x-y) c$,\\[-8pt]
\item If $w_i w_{i+1}=\binom{0\ 1}{1\ 0}$
      then $c:=(xq_iu_{n-i-1}-y)(x-y) c$,\\[-8pt]
\item If $w_i w_{i+1}=\binom{1\ 0}{1\ 0}$
      then $c:=(x-q_{i+1}u_{n-i}y)(x-y) c$,\\[-8pt]
\item If $w_i w_{i+1}=\binom{1\ 1}{1\ 0}$
      then $c:=(xu_{n-i-1}-u_{n-i}y)(x-y) c$,\\[-8pt]
\end{itemize}
where the dot indicates any possible value. Note that if $v_i=0$ and
$v_{i+1}=1$, no factor is added to $c_{IJ}$.

%%%%%%%%%%%%%%%%%%%%%%%%%%%%%%%%%%%%%%%%%%%%%%%%%%%%%%%%%%
\subsubsection{An analogue of the $(1-t)/(1-q)$ transform}

Recall that for commutative symmetric functions, the  $(1-t)/(1-q)$ transform
is defined in terms of the power-sums by
\begin{equation}
p_n\left(\frac{1-t}{1-q}X\right) = \frac{1-t^n}{1-q^n}p_n(X)\,.
\end{equation}
There exist several noncommutative analogues of this transformation. One can
define it on a sequence of generators and require that it be an algebra
morphism. This is the case of the versions chaining internal products
on the right by $\sigma_1((1-t)A)$ and $\sigma_1(A/(1-q))$ (the order
matters).
Taking internal products on the left instead, one obtains linear maps wich are
not algebra morphisms but still lift the commutative transform.

With the specialization $x=1$, $y=t$, $q_i=q^i$, $u_i=1$, $a=b=1$,
one obtains a basis such that for a hook composition $I=(n-k,1^k)$,
the commutative image of $H'_I(U,V)$ becomes the $(1-t)/(1-q)$
transform of the Schur function $s_{n-k,1^k}$.

%%%%%%%%%%%%%%%%%%%%%%%%%%%%%%%%%%%%%%%%%%%%%%%%%%%%%%
\subsubsection{An analogue of the Macdonald $P$-basis}

With the specialization $x=1$, $y=t$, $q_i=q^i$, $u_i=t^i$, $a=b=1$,
one obtains an analogue of the Macdonald $P$-basis, in the sense that for hook
compositions $I=(n-k,1^k)$, the commutative image of $H'_I$ is proportional to
the Macdonald polynomial $P_{n-k,1^k}(q,t;X)$.

For example, the following determinant
\begin{equation}
\left|
\begin {array}{ccccc}
   \left( 1-t \right) h_1          & \left( 1-{t}^{2} \right) {h}_{2}
 & \left( 1-{t}^{3} \right) {h}_{3}& \left( 1-{t}^{4} \right) {h}_{4}
 & \left( 1-{t}^{5} \right) {h}_{5 }\\
\noalign{\medskip}q- 1             & \left( q-t \right) h_1
 & \left( q-{t}^{2} \right) {h}_{2}& \left(q-{t}^ {3} \right) {h}_{3}
 & \left( q-{t}^{4} \right) {h}_{4} \\
\noalign{\medskip}0                &{q}^{2}-1
 & \left( {q}^{2}-t \right) h_1    & \left({q}^{2}-{t}^{2} \right) {h}_{2}
 & \left( {q}^{2}-{t}^{3} \right) {h}_{3} \\
\noalign{\medskip}0                &0
 &1-{t}^{2}                        & \left( 1-{t}^{3} \right) h_1
 & \left(1-{t}^{4} \right) {h}_{2}  \\
\noalign{\medskip}0&0&0&1-t& \left( 1-{t}^{2} \right) h_1
\end {array}
\right|
\end{equation}
is equal to
\begin{equation}
((1-q)(1-q^2))^2 (1-q^5) P_{311}(q,t;X).
\end{equation}
In general, the commutative image of $H'_{n-k,1^k}(U,V)$ is
\begin{equation}
[k]_q!\, [n-k-1]_q!\, (1-q^n)\, P_{n-k,1^k}(q,t;X).
\end{equation}

%%%%%%%%%%%%%%%%%%%%%%%%%%%%%%%%%%%%%%%%%%%%%%%%%%%%%%%%%%%%%%%%%%%%%%%%%%%%%%%
%%%%%%%%%%%%%%%%%%%%%%%%%%%%%%%%%%%%%%%%%%%%%%%%%%%%%%%%%%%%%%%%%%%%%%%%%%%%%%%
%%%%%%%%%%%%%%%%%%%%%%%%%%%%%%%%%%%%%%%%%%%%%%%%%%%%%%%%%%%%%%%%%%%%%%%%%%%%%%%

\bigskip
{\footnotesize\noindent
\sc Institut Gaspard Monge, \\
Universit\'e Paris-Est Marne-la-Vall\'ee,\\
77454 Marne-la-Vall\'ee cedex 2,\\
FRANCE}
\end{document}